\documentclass[11pt]{article}
\usepackage{graphicx}
\usepackage{amsmath}
\usepackage{mathtools}

\usepackage[letterpaper,width=135mm,height=203mm]{geometry}

\newenvironment{Proof}{\textsc{Proof.}}{\hspace{8mm}\QEDmark\smallskip}
\usepackage{latexsym}
\newcommand{\QEDmark}{\mbox{$\Box$}}

\newtheorem{theorem}{Theorem}
\newtheorem{lemma}[theorem]{Lemma}
\newtheorem{conjecture}{Conjecture}

\begin{document}

\begin{center}  {\Large 
\textbf{A Note on One-Hole Domino Tilings of Squares and Rectangles}} \\[2mm]
Seok Hyun Byun and Wayne Goddard \\
School of Mathematical and Statistical Sciences, Clemson University 
\end{center}

\begin{abstract}
We consider the number of domino tilings of an odd-by-odd rectangle that leave one hole. This problem is equivalent to the number of near-perfect matchings of the odd-by-odd rectangular grid. For any particular position of the vacancy on the $(2k+1)\times (2k+1)$ square grid, we show that the number of near-perfect matchings is a multiple of $2^k$, and from this follows a conjecture of Kong that the total number of near-perfect matchings is a multiple of $2^k$. We also determine the parity of the number of near-perfect matchings with a particular vacancy for the rectangle case.
\end{abstract}

\section{Introduction}

Consider a rectangular board. A \textit{domino} is the union of two adjacent unit squares on the board and a \textit{domino tiling} of the board is a collection of dominoes that cover the board without gaps and overlaps. Since every domino consists of two unit squares, if the board has odd area, it cannot have a domino tiling. However, in this case, one can instead think of domino tilings of the board with one vacancy. This paper proves some number-theoretic properties of the total numbers of such tilings.

These problems can be equivalently phrased in terms of graphs and their (near) perfect matchings. Given a graph with an even number of vertices, a \textit{perfect matching} is a collection of edges such that every vertex is contained in exactly one edge. Similarly, given a graph with an odd number of vertices, a \textit{near-perfect matching} is a collection of edges such that all but one vertex is contained in exactly one edge and the one vertex is unmatched. Consider an $m\times n$ rectangular grid. If either $m$ or $n$ is even, then a domino tiling of the board is equivalent to a perfect matching of the grid, and if both $m$ and $n$ are odd, then a domino tiling with one vacancy is equivalent to a near-perfect matching of the grid. Throughout this paper, we phrase our results in terms of grids and their (near) perfect matchings, though one can also state them using the language of domino tilings.

The case when at least one of the side-lengths is even is well-studied. For example, the number of perfect matchings of such a grid is given by an explicit product formula due to Kasteleyn~\cite{kasteleyn} and independently to Temperley and Fisher~\cite{fisher,fishertemperley}. If we restrict ourselves to square grids with even side-length, it is known that the number of perfect matchings exhibits some interesting number-theoretic properties. For instance, Pachter~\cite{pachterDomino} showed that the number of perfect matchings of a square grid of side-length $2k$ is of the form $2^{k}(2\ell +1)^2$, and Cohn~\cite{cohn} further studied the $2$-adic behavior of this number as $k$ varies.

On the other hand, much less is known and studied for the grid with odd side-length. In \cite{kongLattice}, Kong gave a computational approach to the problem of finding the number of near-perfect matchings of the odd square grid, and made some conjectures analogous to the results for the even square grid in \cite{pachterDomino} and \cite{cohn}. Note that Bouttier et al.~\cite{bouttierVacancy} also considered this problem from a more graph-theoretic point of view.

In the first part of the paper, we consider one of the problems posed by Kong in \cite{kongLattice} and give a proof of it. More precisely, for any particular position of the vacancy on the square $(2k+1)\times(2k+1)$ grid, we show that the number of near-perfect matchings is a multiple of $2^k$. From this result, we conclude that the total number of near-perfect matchings is a multiple of $2^k$ and prove one of Kong's conjectures.

In the second half of the paper, we consider a rectangular grid with a fixed vacancy and give some parity and modulo $4$ residue results on the number of near-perfect matchings of it.

\section{The Square Minus a Hole} \label{s:square}

Let $a(2k+1)$ be the number of near-perfect matchings of the $(2k+1) \times (2k+1)$ grid. Kong~\cite{kongLattice} conjectured that $a(2k+1)$ is an odd multiple of $2^k$. We present a proof of this conjecture in Theorems \ref{t:holeyTwos} and \ref{t:main} below.

We need a result about near-perfect matchings with a specific (or fixed) hole. For any particular hole
$h$ of the $(2k+1)\times (2k+1)$ grid, let $a(2k+1; h)$ denote the number of near-perfect matchings with $h$ unmatched, or
equivalently the number of perfect matchings of the subgraph obtaining by deleting $h$.
(Note that in what follows we implicitly assume $h$ has the same color as the corners under the checkerboard coloring, since otherwise there is no 
near-perfect matching, as can be easily seen using a coloring argument.)
This problem has been studied for some particular positions of $h$. For example, when $h$ is on the boundary, by the bijection of Temperley~\cite{temperley}, the value $a(2k+1;h)$ equals the number of spanning trees of the $(k+1) \times (k+1)$ grid, no matter where on the boundary $h$ is. Importantly here, Tenner~\cite{tennerHoley} showed that if $h$ is the center, then $a(2k+1;h)$ is $2^k$ times an odd square (see Corollary~2 in that paper). Also, Bouttier et al.~\cite{bouttierVacancy} generalized Temperley's bijection to show that there is a bijection between (i) near-perfect matchings with a fixed hole, and (ii) spanning webs whose tree component contains the vacancy, together with a choice of orientation for each loop of the spanning web and each loop of the dual web. This bijection will be discussed further and used when we prove Lemma~\ref{l:spanning} in Section~\ref{s:rectangle}.

We now present our first result. In its proof, we use the idea behind the proof of the Matching Factorization Theorem of Ciucu~\cite{ciucuFactor}.

\begin{theorem} \label{t:holeyTwos}
For any particular hole $h$ of the $(2k+1)\times (2k+1)$ grid, $a(2k+1;h)$ is a multiple of $2^k$.
\end{theorem}
\begin{Proof}
By the result of Tenner~\cite{tennerHoley} mentioned above, we know the result is true for the central hole. So assume $h$ is not the center. Fix some diagonal $d$ of the grid that does not contain~$h$.  

Consider a near-perfect matching $M$ with the specified hole $h$. (For an example, see the left picture in Figure~\ref{f:one}.) Let $M'$ be the near-perfect matching that is the mirror image of $M$ around $d$, and let $S$ be the spanning multi-graph with edge set $M \cup M'$. Then, every vertex in $S$ has degree $2$ except for the hole~$h$ and its mirror image $h'$. Thus $S$ consists of a collection of (alternating) cycles (where some cycles might have length $2$) and an (alternating) path $P$ connecting $h'$ to~$h$ (see the right picture in Figure~\ref{f:one}). Further, any cycle of $S$ that contains a vertex of $d$ is symmetric, as is the $h'$--$h$ path. In particular, this means that a cycle that contains a vertex of $d$ contains exactly two of them, and the $h'$--$h$ path contains exactly one vertex of $d$.

\begin{figure}
    \centering \includegraphics[width=1\textwidth]{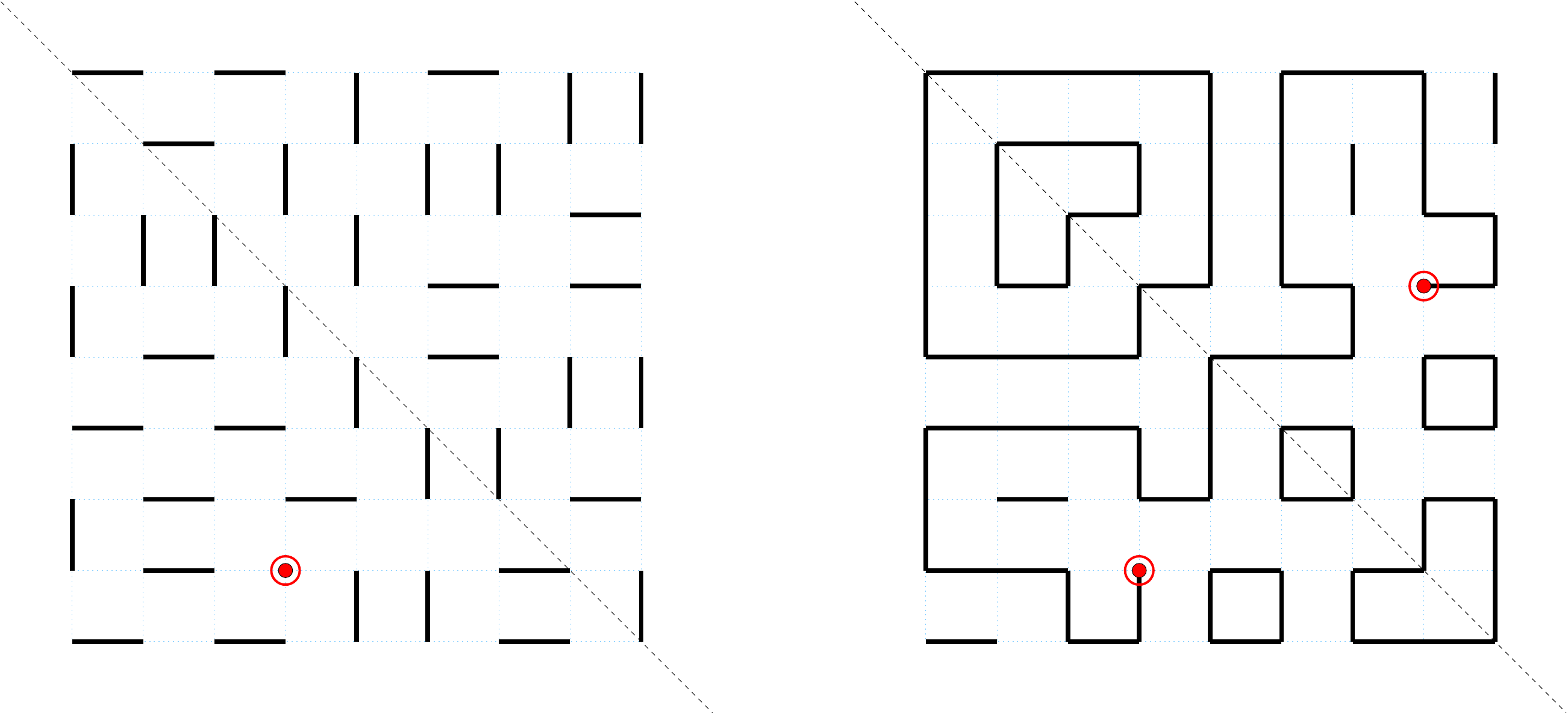}
\caption{The left picture is a near-perfect matching $M$ of $9\times9$ grid with a hole~$h$ (ringed). The right picture shows the union $S=M\cup M'$ ($h$ and $h'$ are ringed). 
}
\label{f:one}
\end{figure}

Let $X$ denote the set of edges in $M$ that are not part of components in $S$ containing vertices of $d$. 
We partition the set of near-perfect matchings with the specified hole $h$ by the pair $(S,X)$. There are $k$ cycles that contain vertices of~$d$, and each cycle can be produced in two ways (depending on which half is~$M$ and which half is~$M'$). So, for each choice of~$(S,X)$, the number of near-perfect matchings is $2^k$. Hence, when one sums up over all choices of $(S,X)$, one obtains that $a(2k+1;h )$ is a multiple of~$2^k$, as required. 
\end{Proof}

Now, for every choice of the hole except the center, there are at least four isomorphic
positions for it. Thus, each orbit collectively contributes a multiple of $2^{k+2}$ to the sum.
So, by the result for the middle hole of Tenner, the overall count is $2^k c_k$ for some odd~$c_k$.
The conjecture of Kong follows:

\begin{theorem} \label{t:main}
For all $k$, the number $a(2k+1)$ of near-perfect matching of a $(2k+1)\times(2k+1)$ square grid  is an odd multiple of $2^k$.
\end{theorem}

Kong~\cite{kongLattice} made some further conjectures on the properties of $c_k$.
For example:

\begin{conjecture}   \label{c:one}
For all $k$,
$a(2k+1) = 2^k c_k$ where $c_k$ is congruent to $1$ modulo $8$. 
\end{conjecture}

We show that Conjecture~\ref{c:one} can be reduced to a much simpler statement.

\begin{theorem}
    If we assume that $a(2k+1;h)$ is a multiple of $2^{k+1}$ for every hole~$h$ that lies on the central row or column but not at the center, then Conjecture~\ref{c:one} holds. 
\end{theorem}

Note that the statement ``$a(2k+1;h)$ is a multiple of $2^{k+1}$ for every hole $h$ that lies on the central row or column but not at the center" has been checked for all integers $k\leq 16$.

\begin{Proof}
    Depending on the position of the hole $h$, we consider the following four cases.
    \begin{itemize}
        \item $h$ is the central hole,
        \item $h$ is on one of the two diagonals but not at the center,
        \item $h$ is on the central row or column but not at the center, and
        \item $h$ is not on any of the symmetry axes.
    \end{itemize}

    If $h$ is the central hole, then by the result of Tenner~\cite{tennerHoley}, $a(2k+1;h)$ is $2^k$ times an odd square. Since an odd square is congruent to $1$ modulo $8$, we have $a(2k+1;h)/2^k\equiv1$ mod $8$ in this case. Thus, it suffices to show that the total number of the near-perfect matchings where the holes is not at the center is a multiple of $2^{k+3}$. If the hole is not on any of the symmetry axes, then there are $8$ isomorphic positions for the hole, and thus, this orbit contributes a multiple of $8\times2^k=2^{k+3}$, as required. On the other hand, if the hole is on the central row or column but not at the center, then there are $4$ isomorphic positions for the hole, and \textbf{by the assumption}, this orbit contributes a multiple of $4\times2^{k+1}=2^{k+3}$. Thus, it suffices to show that the contribution of the second case is also a multiple of $2^{k+3}$.

 Suppose $h$ is a hole on one of the two diagonals but not both. Since there are again $4$ isomorphic positions for this hole, it suffices to show that $a(2k+1;h)$ is a multiple of $2^{k+1}$. Consider a near-perfect matching $M$ with the specific such $h$ (for example, see the left picture in Figure~\ref{f:two}). Let $d'$ be the diagonal that $h$ lies on, and $d$ be the other diagonal. As in the proof of Theorem~\ref{t:holeyTwos}, we reflect $M$ across $d$. Then, we get the spanning multi-graph $S$ that consist of (alternating) cycles and an (alternating) path $P$ connecting $h$ and $h'$ (see the middle picture in Figure~\ref{f:two}). As before, $P$ contains exactly one vertex on $d$ and $S$ contains $k$ cycles that pass through $d$. Thus for any
 specified $S$, the number of near-perfect matchings $M$ of $(2k+1)\times(2k+1)$ that yield $S$ is a multiple of $2^{k}$. 
    
    \begin{figure}
    \centering \includegraphics[width=1\textwidth]{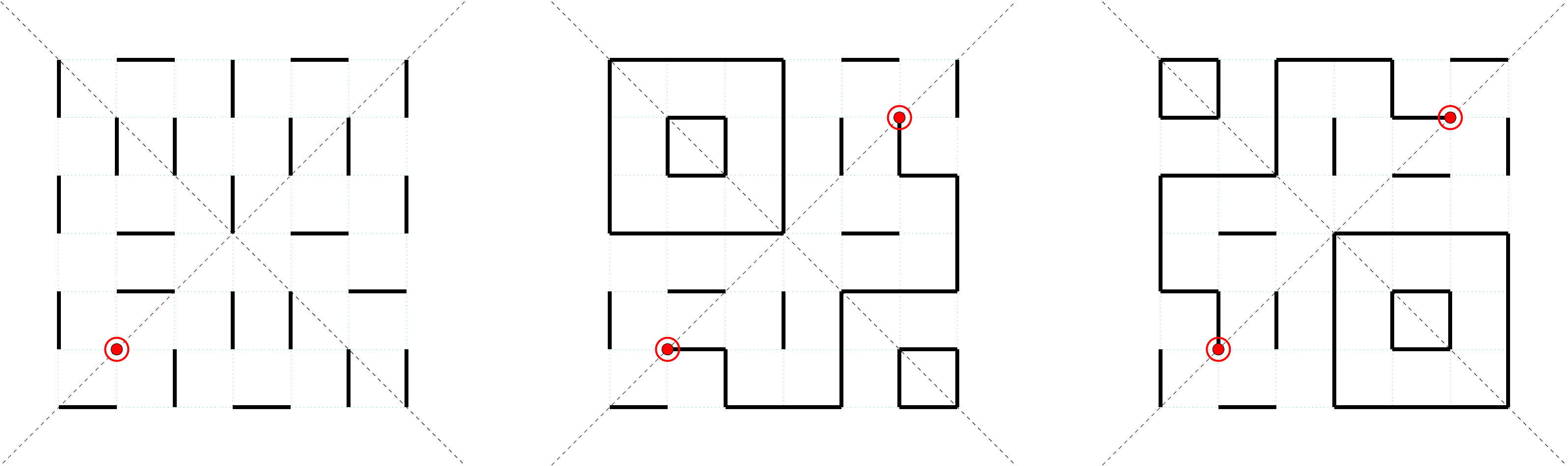}
\caption{The left picture shows a matching $M$ of $7\times7$ grid with hole $h$ (ringed). The spanning multi-graph $S$ with edge set $M\cup M'$ is shown in the middle picture. The mirror image $S'$ is shown in the right picture. 
}
\label{f:two}
\end{figure}

 Now,  observe that the configuration $S$ cannot be symmetric across $d'$, since vertex $h$ has only one incident edge. Let $S'$ be the mirror image of $S$ across $d'$ (see the right picture in Figure~\ref{f:two}). Then, by symmetry, $S'$ consists of the mirror images of the cycles in $S$ and a path $P'$ connecting $h$ and $h'$ that is a mirror image of the path $P$. Thus, the number of near-perfect matchings that yield $S'$ is equal to the number that yield $S$. Since the number $a(2k+1,h)$ can be written as a sum over $S$ of the number of matchings $M$ such that $M\cup M'$ is equal to $S$ or $S'$, the number $a(2k+1,h)$ is a multiple of $2^{k+1}$. This completes the proof.
\end{Proof}

%\newpage

\section{A Rectangle minus a Hole} \label{s:rectangle}

The question for a rectangular grid with a hole is less explored. Here the total number of near-perfect matchings can be odd. So an obvious question is: when is it odd? We address this by again looking at the question for specific holes. For odd $r$ and $c$, let $a(r,c)$ be the number of near-perfect matchings of the $r\times c$ grid; and  for any particular hole $h$, let $a(r,c; h)$ denote the number of near-perfect matchings that leave that specified hole unmatched. Note that due to the reduced symmetry of the rectangle, most of the arguments used in the previous section are not applicable in this case.

Say we have a $r\times c$ grid with $r$ and $c$ odd imbued with a checkerboard coloring. We say a vertex is \textit{white} if it is of the majority coloring; otherwise, it is \textit{black}. Further, a white vertex is \textit{odd} or \textit{even} based on the parity of its row or column (so that white vertices on the boundary are odd).

\begin{theorem} \label{t:rectangleOdd}
For a white hole $h$, the number of near-perfect matchings $a(r,c;h)$ is odd if and only if the following three conditions all hold:\\
(i) At least one of $r$ or $c$ is congruent to $1$ mod $4$;  \\
(ii) $\gcd(r+1,c+1)=2$; and \\
(iii) the hole $h$ is odd.
\end{theorem}

Our proof of Theorem~\ref{t:rectangleOdd} consists of three parts. The first part is to show that if neither condition (i) nor condition (ii) holds,
then $a(r,c;h)$ is always even. We need an old theorem of Lov\'asz (see \cite{Lovasz}):

\begin{theorem} \label{t:Lovasz}
The number of perfect matchings in any graph is even if and only if there
is a nonempty subset $S$ of the vertices such that every vertex has an even number of neighbors 
in $S$.
\end{theorem}

\begin{figure}
    \centering \includegraphics{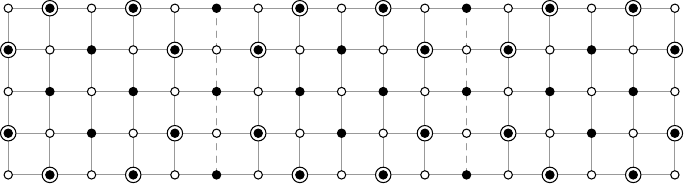}
\caption{A subset of the $5\times 17$ grid (ringed) where each vertex has an even number of neighbors in the subset.} 
\label{f:5by17}
\end{figure}

\begin{lemma} \label{l:evenIandII}
(a) If both $r$ and $c$ are congruent to $3$ mod $4$, then $a(r,c;h)$ is even.\\
(b) If both $r$ and $c$ are odd and $\gcd(r+1,c+1)$ is divisible by some odd $f>1$, then $a(r,c;h)$ is even.
\end{lemma}
\begin{Proof}
(a) We assign $(x,y)$-coordinates to vertices in a $r$-by-$c$ rectangle such that the vertex at the bottom left corner has a coordinate $(1,1)$. Let $S$ be all vertices with coordinates $(i+1, j)$, $(i-1,j)$, $(i,j-1)$ and $(i,j+1)$ where $i,j \equiv 2 \mbox{\,mod\,} 4$. 
Figure~\ref{f:7by11} gives an example of $S$ with $r=7$ and $c=11$.
It can readily be checked that $S$ has the desired property in Theorem~\ref{t:Lovasz}; that is, every vertex has an even number of neighbors in $S$. Now, note that all vertices in $S$ are black (their coordinates sum to odd). Thus, the set $S$ also satisfies the desired property after removing any collection of white vertices from the grid. In particular, $a(r,c;h)$ is even for any white hole $h$. 
   \begin{figure}
    \centering \includegraphics{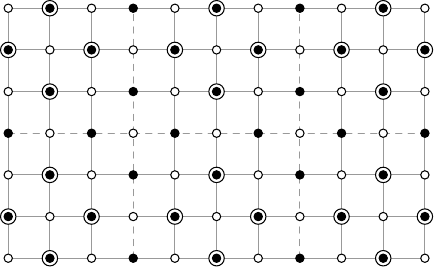}
\caption{A subset of the $7\times 11$ grid (ringed) where each vertex in the grid has an even number of neighbors in the subset.} 
\label{f:7by11}
\end{figure}

(b) Note that both $r+1$ and $c+1$ are multiples of $2f$. First partition the rectangle into $(2f-1) \times (2f-1)$ subgraphs with one empty row and one empty column between consecutive subgraphs. Then form $S$ by, in each subgraph, taking the $4(f-1)$ vertices whose Manhattan distance from the center is $f$. 
Figure~\ref{f:5by17} is an example with $f=3$: namely $r=5$ and $c=17$.
Again, $S$ satisfies the condition of Theorem~\ref{t:Lovasz} and all vertices in $S$ are black. Hence, again the graph with any white vertex removed has an even number of perfect matchings. 
\end{Proof}

We continue the proof of Theorem~\ref{t:rectangleOdd}.

\begin{figure}
    \centering \includegraphics[width=0.75\textwidth]{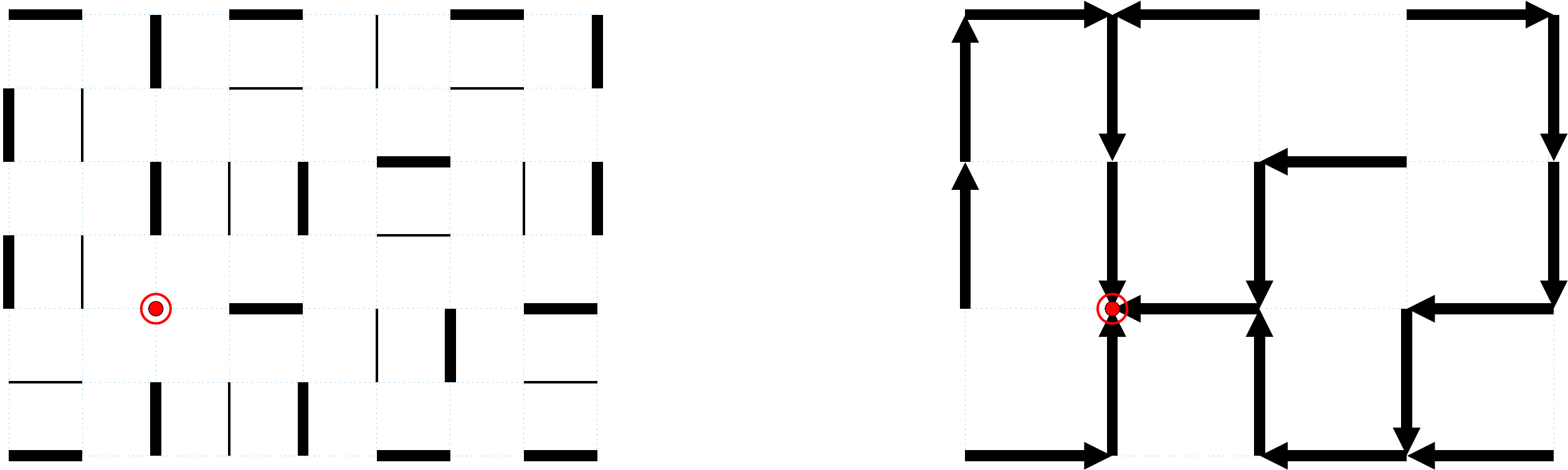}
\caption{
The left picture shows a near-perfect matching $M$ of  the $7\times9$ grid with a vacancy $h$ (ringed) an odd white vertex. The corresponding directed spanning web $S_o$ is shown in the right picture.}
\label{f:webTwo}
\end{figure}

\begin{lemma} \label{l:spanning}
For all odd $r$ and $c$, \\
(a) If the white hole $h$ is even, then  $a(r,c;h)$ is even. \\
(b) If the white hole $h$ is odd, then $a(r,c;h)$ has the same parity
regardless of where the hole actually is.
\end{lemma}
\begin{Proof}
We use the ideas behind spanning webs (see, for example, Bouttier et al.~\cite{bouttierVacancy}). Specifically, associated with the $r \times c$ grid $R$, we define a graph $H_o$ whose vertices are the odd white vertices of $R$, where two vertices are adjacent if they are two apart in the same row or column of $R$. As noted for example in~\cite{bouttierVacancy}, given a near-perfect matching $M$ of $R$, one can associate a directed spanning subgraph $S_o$ of $H_o$, called a \textit{spanning web}. Specifically, for each odd white vertex~$w_1$ consider the edge $e$ in $M$ (if any) that covers $w_1$. 
(See the left picture in Figure~\ref{f:webOne}, where the thickened edges are the ones covering odd white vertices.)
Say $e$ also covers black vertex~$b$. Then let $w_2$ be the other neighbor of $b$ in the same row/column (exists since all white vertices on the outside of $R$ are odd); and associate this edge $e$ in $M$ with the arc $w_1 w_2$ directed from $w_1$ to~$w_2$. This produces a directed subgraph $S_o$ of $H_o$ where every covered white vertex has an out-degree one, and every uncovered white vertex has an out-degree zero.

\begin{figure}
    \centering \includegraphics[width=1\textwidth]{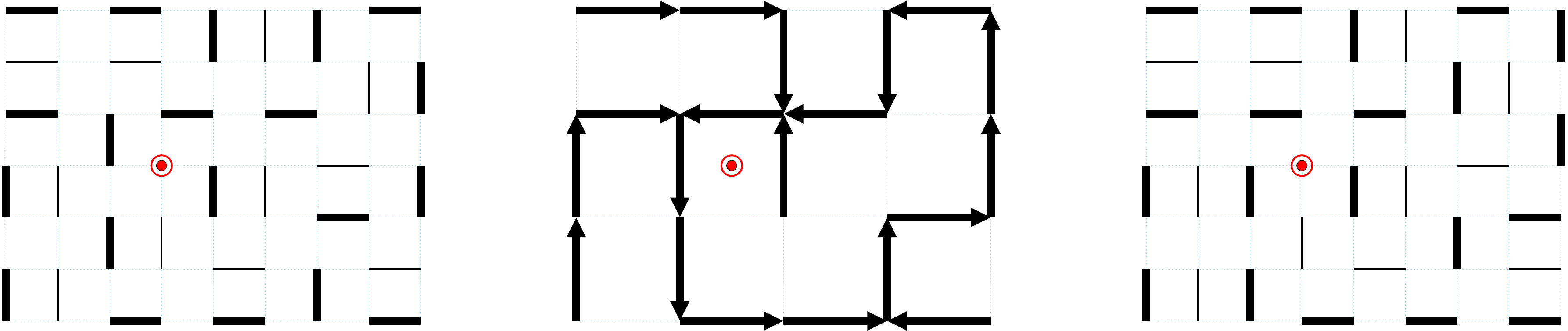}
\caption{The left picture show a near-perfect matching $M$ of a $7\times9$ grid with the vacancy $h$ (ringed) an even white vertex. The middle picture shows the corresponding spanning web $S_o$ (in a $4\times5$ grid). The right picture shows the matching $M'$ obtained from $M$ by reversing the cycle in the middle picture.}
\label{f:webOne}
\end{figure}

Assume first that the hole $h$ is even. Consider a near-perfect matching $M$ of~$R$. Then every vertex of~$S_o$ has out-degree exactly one, and thus the corresponding spanning web necessarily contains a directed cycle. This cycle corresponds to an alternating cycle in $R$. 
If one reverses the cycle in $S_o$, then one obtains the ``complementary'' near-perfect matching $M'$ (see the right picture in Figure~\ref{f:webOne}). 
A parity argument (for example noted in~\cite{bouttierVacancy}) shows that all cycles surround~$h$. 
Thus one can pair off the near-perfect matchings, say by their innermost cycle. 
Hence the total number of near-perfect matchings is even.

Assume second that the hole $h$ is odd. Then every vertex of~$S_o$ has out-degree one, except for $h$ itself, which has out-degree zero. The directed subgraph $S_o$ might or might not have a cycle. But again one can pair off the near-perfect matchings for which there is a directed cycle, based on reversing say the innermost cycle. So to determine the parity, one can restrict to the near-perfect matchings where the spanning web has no cycle, and so is a spanning tree of $H_o$ rooted at~$h$. (See Figure~\ref{f:webTwo}.) As noted in~\cite{bouttierVacancy} (or due to Temperley's bijection~\cite{temperley}), there is a unique near-perfect matching for each such $S_o$. So the total number of near-perfect matchings has the same parity as the number of rooted spanning trees; but the number of rooted spanning trees with a fixed root is just the number of spanning trees, and is thus the same for all choices of $h$.
\end{Proof}

We next determine the parity in one particular case:

\begin{lemma} \label{l:fold}
Assume $c$ is congruent to $1$ mod $4$ and $\gcd(r+1,c+1)=2$. For the hole $h^*$ in the middle of the first row, it holds that $a(r,c;h^*)$ is odd.
\end{lemma}
\begin{Proof}
We count the number of near-perfect matchings. Those that are not symmetric with respect to the middle column can be paired off. So $a(r,c;h^*)$ has the same parity as the number of near-perfect matchings $M$ that are symmetric with respect to the middle column. 

This means that the middle column except for $h^*$ is filled by vertical edges of $M$. Hence this number equals the number of perfect matchings of the $c' = (c-1)/2$ by $r$ subrectangle that results. Note that $c'$ is even and so $\gcd( c'+1, r+1)$ is odd. At the same time, $\gcd( c'+1, r+1)$ divides $ \gcd( 2(c'+1), r+1) ) = \gcd ( c+1, r+1 ) = 2$. That is, $\gcd( c'+1, r+1)=1$. So by Proposition 5.10 of Barkley and Liu~\cite{barkleyChannel},  the number of perfect matchings of such a rectangle is odd.
\end{Proof}

By Lemma~\ref{l:evenIandII} and \ref{l:spanning}(a), the number of near perfect matchings is even if any of the three conditions in the theorem fail. By Lemma~\ref{l:spanning}(b), in the remaining case the parity of the number of near-perfect matching is independent of the exact placement of the hole, and by Lemma~\ref{l:fold} this is odd. This completes the proof of Theorem~\ref{t:rectangleOdd}.

Since the total number of near-perfect matchings is odd if and only if the first two conditions of Theorem~\ref{t:rectangleOdd} hold and 
the number $(r+1)(c+1)/4$ of odd white vertices is odd,
as a consequence we get:

\begin{theorem}
For odd $r,c$, the number $a(r,c)$  of near-perfect matchings of the $r\times c$ grid  is odd if and only if both $r$ and $c$ are congruent to $1$ mod $4$,
and $\gcd(r+1,c+1)=2$.
\end{theorem}

A natural next question is to consider counts modulo $4$. 
For example, the proof of Lemma~\ref{l:spanning}(b) can be extended to show that the number of near perfect matchings for an odd white hole is the same modulo $4$ for all odd white holes. This uses the fact that in this case there is also guaranteed to be a cycle in $S_e$, the graph analogous to $S_o$ except that its vertices are the even white vertices. 
We omit the details.

We conclude with one example conjecture:

\begin{conjecture}
Except sometimes when $r$ and $c$ are both congruent to $3$ modulo~$4$, the count $a(r,c;h)$ is a multiple of $4$ for any even white hole $h$.
\end{conjecture}

\section*{Acknowledgments}

The authors would like to thank the referee for the valuable comments.

\end{document}